\documentclass[12pt,vatola]{article}

\usepackage{amsfonts}

\usepackage{mathrsfs}

\usepackage{graphicx}

\def\c{\centerline}

\def\re#1{\par\hangindent\parindent\indent\llap{#1\enspace}\ignorespaces}

\def\no{\noindent}

\begin{document}

\c{\bf\large A Generalization of}\vskip 3mm

\c{\bf\large Seifert-Van Kampen Theorem for Fundamental Groups}
\vskip 8mm

\c{Linfan MAO}\vskip 5mm

\c{\scriptsize (Chinese Academy of Mathematics and System Science,
Beijing 100080, P.R.China)}

\c{\scriptsize E-mail: maolinfan@163.com}

\vskip 12mm

\begin{minipage}{130mm}

\no{\bf Abstract}: {\small As we known, the {\it Seifert-Van Kampen
theorem} handles fundamental groups of those topological spaces
$X=U\cup V$ for open subsets $U,\ V\subset X$ such that $U\cap V$ is
arcwise connected. In this paper, this theorem is generalized to
such a case of maybe not arcwise-connected, i.e., there are $C_1$,
$C_2$,$\cdots,\ C_m$ arcwise-connected components in $U\cap V$ for
an integer $m\geq 1$, which enables one to find fundamental groups
of combinatorial spaces by that of spaces with theirs underlying
topological graphs, particularly, that of compact manifolds by their
underlying graphs of charts.}\vskip 2mm

\no{\bf Key Words}: {\small Fundamental group, Seifert-Van Kampen
theorem, topological space, combinatorial manifold, topological
graph.}\vskip 2mm

\no{\bf AMS(2010)}: {\small 51H20.}

\end{minipage}

\vskip 12mm

\no{\bf \S $1.$ \ Introduction}

\vskip 5mm

\no All spaces $X$ considered in this paper are arcwise-connected,
graphs are connected topological graph, maybe with loops or multiple
edges and interior of an arc $a: (0,1)\rightarrow X$ is opened. For
terminologies and notations not defined here, we follow the
reference [1]-[3] for topology and [4]-[5] for topological graphs.

Let $X$ be a topological space. A {\it fundamental group}
$\pi_1(X,x_0)$ of $X$ based at a point $x_0\in X$ is formed by
homotopy arc classes in $X$ based at $x_0\in X$. For an
arcwise-connected space $X$, it is known that $\pi_1(X,x_0)$ is
independent on the base point $x_0$, that is, for $\forall x_0,
y_0\in X$,
$$\pi_1(X,x_0)\cong\pi_1(X,y_0).$$

Find the fundamental group of a space $X$ is a difficult task in
general. Until today, the basic tool is still the {\it Seifert-Van
Kampen theorem} following.\vskip 4mm

\no{\bf Theorem $1.1$}(Seifert and Van-Kampen) \ {\it Let $X=U\cup
V$ with $U,\ V$ open subsets and let $X,\ U,\ V$, $U\cap V$ be
non-empty arcwise-connected with $x_0\in U\cap V$ and $H$ a group.
If there are homomorphisms}
$$\phi_1:\pi_1(U,x_0)\rightarrow H \ \ {and} \ \ \phi_2:\pi_1(V,x_0)\rightarrow H$$

\no{\it and}

\c{\includegraphics[bb=10 10 180 140]{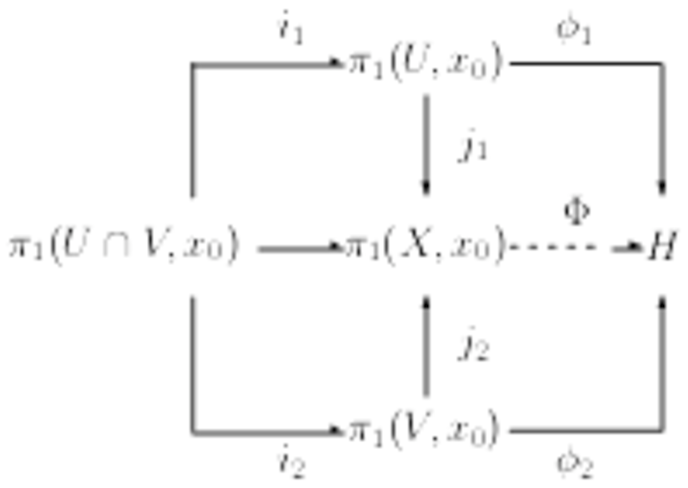}}\vskip 3mm

\no{\it with $\phi_1\cdot i_1=\phi_2\cdot i_2$, where
$i_1:\pi_1(U\cap V,x_0)\rightarrow\pi_1(U,x_0)$, $i_2:\pi_1(U\cap
V,x_0)\rightarrow\pi_1(V,x_0)$,
$j_1:\pi_1(U,x_0)\rightarrow\pi_1(X,x_0)$ and
$j_2:\pi_1(V,x_0)\rightarrow\pi_1(X,x_0)$ are homomorphisms induced
by inclusion mappings, then there exists a unique homomorphism
$\Phi:\ \pi_1(X,x_0)\rightarrow H$ such that $\Phi\cdot j_1=\phi_1$
and $\Phi\cdot j_2=\phi_2$.} \vskip 2mm

Applying Theorem $1.1$, it is easily to determine the fundamental
group of such spaces $X=U\cup V$ with $U\cap V$ an arcwise-connected
following.

\vskip 3mm

\no{\bf Theorem $1.2$}(Seifert and Van-Kampen theorem, classical
version) \ {\it Let spaces $X,U,V$ and $x_0$ be in Theorem $1.1$.
If}
$$j: \pi_1(U,x_0)*\pi_1(V,x_0)\rightarrow\pi_1(X,x_0)$$

\no{\it is an extension homomorphism of $j_1$ and $j_2$, then $j$ is
an epimorphism with kernel {\rm Ker$j$} generated by
$i_1^{-1}(g)i_2(g),\ g\in\pi_1(U\cap V,x_0)$, i.e.,}
$$\pi_1(X,x_0)\cong\frac{\pi_1(U,x_0)*\pi_1(V,x_0)}
{\left[i_1^{-1}(g)\cdot i_2(g)|\ g\in\pi_1(U\cap V,x_0)\right]},$$
where $\left[A\right]$ denotes the minimal normal subgroup of a
group $\mathscr{G}$ included $A\subset\mathscr{G}$.

\vskip 3mm

Now we use the following convention.\vskip 3mm

\no{\bf Convention $1.3$} \ {\it Assume that \vskip 2mm

{\rm($1$)} $X$ is an arcwise-connected spaces, $x_0\in X$;\vskip 1mm

{\rm($2$)} $\{U_{\lambda}: \lambda\in\Lambda\}$ is a covering of $X$
by arcwise-connected open sets such that $x_0\in U_{\lambda}$ for
$\forall\lambda\in\Lambda$;\vskip 1mm

{\rm($3$)} For any two indices $\lambda_1,\lambda_2\in\Lambda$ there
exists an index $\lambda\in\Lambda$ such that $U_{\lambda_1}\cap
U_{\lambda_2}=U_{\lambda}$}\vskip 2mm

If $U_{\lambda}\subset U_{\mu}\subset X$, then the notation
$$\phi_{\lambda\mu}:\pi_1(U_{\lambda},x_0)\rightarrow\pi_1(U_{\mu},x_0) \ \ {\rm and} \ \
\phi_{\lambda}:\pi_1(U_{\lambda},x_0)\rightarrow\pi_1(X,x_0)$$
denote homomorphisms induced by the inclusion mapping
$U_{\lambda}\rightarrow U_{\mu}$ and $U_{\lambda}\rightarrow X$,
respectively. It should be noted that the Seifert-Van Kampen theorem
has been generalized in [2] under Convention $1.3$ by any number of
open subsets instead of just by two open subsets following.

\vskip 4mm

\no{\bf Theorem $1.4$}([2]) \ {\it Let $X, U_{\lambda},\
\lambda\in\Lambda$ be arcwise-connected space with Convention $1.3$
satisfies the following universal mappping condition: Let $H$ be a
group and let $\rho_{\lambda}:\pi_1(U_{\lambda},x_0)\rightarrow H$
be any collection of homomorphisms defined for all
$\lambda\in\Lambda$ such that the following diagram is commutative
for $U_{\lambda}\subset U_{\mu}$:}

\c{\includegraphics[bb=10 10 180 80]{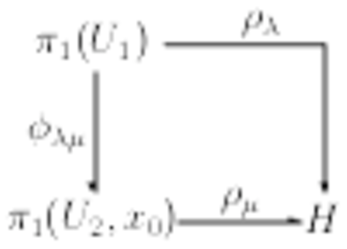}}\vskip 3mm

\no{\it Then there exists a unique homomorphism
$\Phi:\pi_1(X,x_0)\rightarrow H$ such that for any
$\lambda\in\Lambda$ the following diagram is commutative:}

\c{\includegraphics[bb=10 10 180 80]{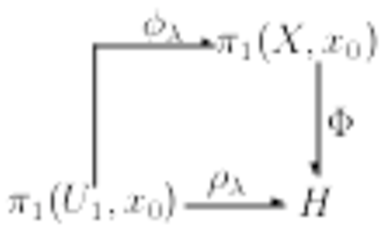}}\vskip 3mm

\no{\it Moreover, this universal mapping condition characterizes
$\pi_1(X,x_0)$ up to a unique isomorphism.}

\vskip 3mm

Theorem $1.4$ is useful for determining the fundamental groups of
CW-complexes, particularly, the adjunction of $n$-dimensional cells
to a space for $n\geq 2$. Notice that the essence in Theorems $1.2$
and $1.4$ is that $\cap_{\lambda\in\Lambda}U_{\lambda}$ is
arcwise-connected, which limits the application of such kind of
results. The main object of this paper is to generalize the
Seifert-Van Kampen theorem to such an intersection maybe non-arcwise
connected, i.e., there are $C_1$, $C_2$,$\cdots,\ C_m$
arcwise-connected components in $U\cap V$ for an integer $m\geq 1$.
This enables one to determine the fundamental group of topological
spaces, particularly, combinatorial manifolds introduced in [6]-[8]
following which is a special case of Smarandache multi-space
([9]-[10]).

\vskip 4mm

\no{\bf Definition $1.4$} \ {\it A combinatorial Euclidean space
$\mathscr{E}_G(n_{\nu};\nu\in\Lambda)$ underlying a connected graph
$G$ is a topological spaces consisting of ${\bf R}^{n_{\nu}}$,
$\nu\in\Lambda$ for an index set $\Lambda$ such that\vskip 3mm

$V(G)=\{{\bf R}^{n_{\nu}}|\nu\in\Lambda\}$;\vskip 2mm

$E(G)=\{\ ({\bf R}^{n_{\mu}},{\bf R}^{n_{\nu}}) | \ {\bf
R}^{n_{\mu}}\cap{\bf R}^{n_{\nu}}\not=\emptyset,
\mu,\nu\in\Lambda\}$.}\vskip 3mm

A {\it combinatorial fan-space} $\widetilde{\bf
R}(n_{\nu};\nu\in\Lambda)$ is a combinatorial Euclidean space
$\mathscr{E}_{K_{|\Lambda|}}(n_{\nu};\nu\in\Lambda)$ of ${\bf
R}^{n_{\nu}},\ \nu\in\Lambda$ such that for any integers
$\mu,\nu\in\Lambda,\ \mu\not= \nu$,
$${\bf R}^{n_{\mu}}\bigcap{\bf R}^{n_{\nu}}= \bigcap\limits_{\lambda\in\Lambda}{\bf
R}^{n_{\lambda}},$$ which enables us to generalize the conception of
manifold to combinatorial manifold, also a locally combinatorial
Euclidean space following. \vskip 4mm

\no{\bf Definition $1.5$} \ {\it For a given integer sequence $0<
n_1<n_2<\cdots< n_m$, $m\geq 1$, a topological combinatorial
manifold $\widetilde{M}$ is a {\it Hausdorff space} such that for
any point $p\in \widetilde{M}$, there is a local chart
$(U_p,\varphi_p)$ of $p$, i.e., an open neighborhood $U_p$ of $p$ in
$\widetilde{M}$ and a homoeomorphism $\varphi_p:
U_p\rightarrow\widetilde{\bf R}(n_1(p),
n_2(p),\cdots,n_{s(p)}(p))=\bigcup\limits_{i=1}^{s(p)}{\bf
R}^{n_i(p)}$ with $\{n_1(p),
n_2(p),\cdots,n_{s(p)}(p)\}\subseteq\{n_1,n_2,\cdots,n_m\}$ and
$\bigcup\limits_{p\in\widetilde{M}}\{n_1(p),
n_2(p),\cdots,n_{s(p)}(p)\}=\{n_1,n_2,\cdots,n_m\}$, denoted by
$\widetilde{M}(n_1,n_2,\cdots,n_m)$ or $\widetilde{M}$ on the
context and
$$\widetilde{{\mathcal A}}=\{(U_p,\varphi_p)|
p\in\widetilde{M}(n_1,n_2,\cdots,n_m))\}$$
an atlas on
$\widetilde{M}(n_1,n_2,\cdots,n_m)$.

A topological combinatorial manifold
$\widetilde{M}(n_1,n_2,\cdots,n_m)$ is finite if it is just combined
by finite manifolds without one manifold contained in the union of
others.}\vskip 3mm

If these manifolds $M_i, \ 1\leq i\leq m$ in
$\widetilde{M}(n_1,n_2,\cdots,n_m)$ are Euclidean spaces ${\bf
R}^{n_i}, \ 1\leq i\leq m$, then $\widetilde{M}(n_1,n_2,\cdots,n_m)$
is nothing but a combinatorial Euclidean space
$\mathscr{E}_{G}(n_{\nu};\nu\in\Lambda)$ with
$\Lambda=\{1,2,\cdots,m\}$. Furthermore, If $m=1$ and $n_1=n$, or
$n_{\nu}=n$ for $\nu\in\Lambda$, then
$\widetilde{M}(n_1,n_2,\cdots,n_m)$ or
$\mathscr{E}_{G}(n_{\nu};\nu\in\Lambda)$ is exactly a manifold
$M^{n}$ by definition.

\vskip 10mm

\no{\bf \S $2.$ \ Topological Space Attached Graphs}

\vskip 4mm

\no A topological graph $G$ is itself a topological space formally
defined as follows.

\vskip 4mm

\no{\bf Definition $2.1$} \ {\it A topological graph $G$ is a pair
$(S,S^0)$ of a Hausdorff space $S$ with its a subset $S^0$ such
that}\vskip 2mm

($1$) {\it $S^0$ is discrete, closed subspaces of $S$;}\vskip 1mm

($2$) {\it $S-S^0$ is a disjoint union of open subsets
$e_1,e_2,\cdots,e_m$, each of which is homeomorphic to an open
interval $(0,1)$;}\vskip 1mm

($3$) {\it the boundary $\overline{e}_i-e_i$ of $e_i$ consists of
one or two points. If $\overline{e}_i-e_i$ consists of two points,
then $(\overline{e}_i,e_i)$ is homeomorphic to the pair
$([0,1],(0,1))$; if $\overline{e}_i-e_i$ consists of one point, then
$(\overline{e}_i,e_i)$ is homeomorphic to the pair
$(S^1,S^1-\{1\})$;}\vskip 1mm

($4$) {\it a subset $A\subset G$ is open if and only if
$A\cap\overline{e}_i$ is open for $1\leq i\leq m$.}

\vskip 3mm

\c{\includegraphics[bb=10 10 180 160]{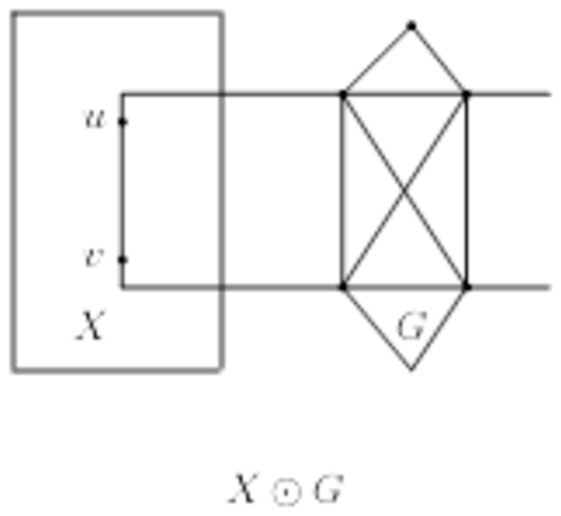}}\vskip 3mm

\c{\bf Fig.$2.1$}\vskip 2mm

Notice that a topological graph maybe with semi-edges, i.e., those
edges $e^+\in E(G)$ with $e^+: [0,1) \ {\rm or} \ (0,1]\rightarrow
S$. A topological space $X$ attached with a graph $G$ is such a
space $X\odot G$ such that
$$X\bigcap G\not=\emptyset, \ \ G\not\subset X$$
and there are semi-edges $e^+\in (X\bigcap G)\setminus G$, $e^+\in
G\setminus X$. An example for $X\odot G$ can be found in Fig.$2.1$.
In this section, we characterize the fundamental groups of such
topological spaces attached with graphs.

\vskip 4mm

\no{\bf Theorem $2.2$} \ {\it Let $X$ be arc-connected space, $G$ a
graph and $H$ the subgraph $X\cap G$ in $X\odot G$. Then for $x_0\in
X\cap G$,}
$$\pi_1(X\odot G,x_0)\cong\frac{\pi_1(X,x_0)*\pi_1(G,x_0)}
{\left[i_1^{-1}(\alpha_{e_{\lambda}})i_2(\alpha_{e_{\lambda}})|\
e_{\lambda}\in E(H)\setminus T_{span})\right]},$$ {\it where $i_1:
\pi_1(H,x_0)\rightarrow X$, $i_2:\pi_1(H,x_0)\rightarrow G$ are
homomorphisms induced by inclusion mappings, $T_{span}$ is a
spanning tree in $H$,
$\alpha_{\lambda}=A_{\lambda}e_{\lambda}B_{\lambda}$ is a loop
associated with an edge $e_{\lambda}=a_{\lambda}b_{\lambda}\in
H\setminus T_{span}$,  $x_0\in G$ and $A_{\lambda}$, $B_{\lambda}$
are unique paths from $x_0$ to $a_{\lambda}$ or from $b_{\lambda}$
to $x_0$ in $T_{span}$.}

\vskip 3mm

{\it Proof} \ This result is an immediately conclusion of
Seifert-Van Kampen theorem. Let $U=X$ and $V=G$. Then $X\odot
G=X\cup G$ and $X\cap G=H$. By definition, there are both semi-edges
in $G$ and $H$. Whence, they are opened. Applying the Seifert-Van
Kampen theorem, we get that
$$\pi_1(X\odot G,x_0)\cong\frac{\pi_1(X,x_0)*\pi_1(G,x_0)}
{\left[i_1^{-1}(g)i_2(g)|\ g\in\pi_1(X\cap G,x_0)\right]},$$

Notice that the fundamental group of a graph $H$ is completely
determined by those of its cycles ([2]), i.e.,\vskip 3mm

\c{$\pi_1(H,x_0)=\left<\alpha_{\lambda}|e_{\lambda}\in E(H)\setminus
T_{span}\right>,$}\vskip 2mm

\no where $T_{span}$ is a spanning tree in $H$,
$\alpha_{\lambda}=A_{\lambda}e_{\lambda}B_{\lambda}$ is a loop
associated with an edge $e_{\lambda}=a_{\lambda}b_{\lambda}\in
H\setminus T_{span}$,  $x_0\in G$ and $A_{\lambda}$, $B_{\lambda}$
are unique paths from $x_0$ to $a_{\lambda}$ or from $b_{\lambda}$
to $x_0$ in $T_{span}$. We finally get the following conclusion,
$$\hskip 20mm\pi_1(X\odot G,x_0)\cong\frac{\pi_1(X,x_0)*\pi_1(G,x_0)}
{\left[i_1^{-1}(\alpha_{e_{\lambda}})i_2(\alpha_{e_{\lambda}})|\
e_{\lambda}\in E(H)\setminus T_{span})\right]}\hskip 20mm\Box$$

\vskip 3mm

\no{\bf Corollary $2.3$} \ {\it Let $X$ be arc-connected space, $G$
a graph. If $X\cap G$ in $X\odot G$ is a tree, then
$$\pi_1(X\odot G,x_0)\cong\pi_1(X,x_0)*\pi_1(G,x_0).$$
Particularly, if $G$ is graphs shown in Fig.$2.2$ following}\vskip
3mm

\c{\includegraphics[bb=10 10 200 100]{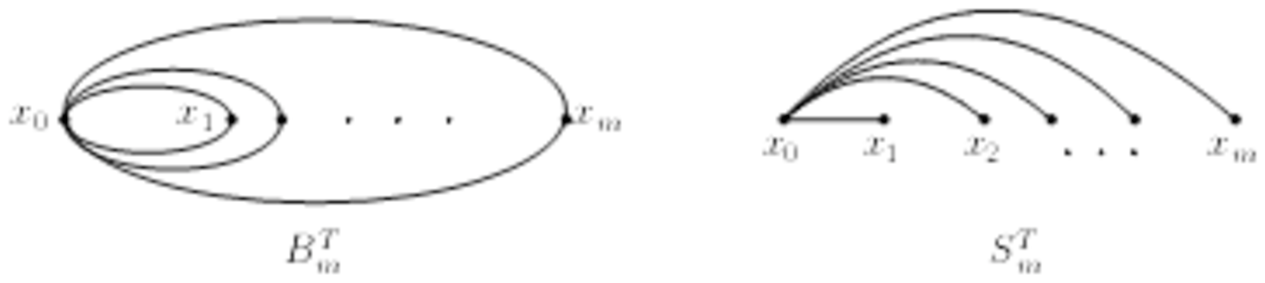}}\vskip 3mm

\c{\bf Fig.$2.2$}\vskip 2mm

\no{\it and $X\cap G= K_{1,m}$, Then}
$$\pi_1(X\odot B_m^T,x_0)\cong\pi_1(X,x_0)*\left<L_i|1\leq i\leq m\right>,$$
{\it where $L_i$ is the loop of parallel edges $(x_0,x_i)$ in
$B_m^T$ for $1\leq i\leq m-1$ and}
$$\pi_1(X\odot S_m^T,x_0)\cong\pi_1(X,x_0).$$

\vskip 3mm

\no{\bf Theorem $2.4$} \ {\it Let $\mathscr{X}_m\odot G$ be a
topological space consisting of $m$ arcwise-connected spaces $X_1,
X_2,\cdots, X_m$, $X_i\cap X_j=\emptyset$ for $1\leq i,j\leq m$
attached with a graph $G$, $V(G)=\{x_0, x_1,\cdots, x_{l-1}\}$,
$m\leq l$ such that $X_i\cap G=\{x_i\}$ for $0\leq i\leq l-1$. Then}
\begin{eqnarray*}
\pi_1(\mathscr{X}_m\odot
G,x_0)&\cong&\left(\prod\limits_{i=1}^m\pi_1(X_i^*,x_0)\right)*\pi_1(G,x_0)\\
&\cong&\left(\prod\limits_{i=1}^m\pi_1(X_i,x_i)\right)*\pi_1(G,x_0),
\end{eqnarray*}
{\it where $X_i^*=X_i\bigcup (x_0,x_i)$ with $X_i\cap
(x_0,x_i)=\{x_i\}$ for $(x_0,x_i)\in E(G)$, integers $1\leq i\leq
m$.}\vskip 3mm

{\it Proof} \ The proof is by induction on $m$. If $m=1$, the result
is hold by Corollary $2.3$.

Now assume the result on $\mathscr{X}_m\odot G$ is hold for $m\leq
k< l-1$. Consider $m=k+1\leq l$. Let $U=\mathscr{X}_k\odot G$ and
$V=X_{k+1}$. Then we know that $\mathscr{X}_{k+1}\odot G=U\cup V$
and $U\cap V=\{x_{k+1}\}$.

Applying the Seifert-Van Kampen theorem, we find that
\begin{eqnarray*}
\pi_1(\mathscr{X}_{k+1}\odot
G,x_{k+1})&\cong&\frac{\pi_1(U,x_{k+1})*\pi_1(V,x_{k+1})}
{\left[i_1^{-1}(g)i_2(g)|\ g\in\pi_1(U\cap V,x_{k+1})\right]}\\
&\cong&\frac{\pi_1(\mathscr{X}_k\odot G,x_0)*\pi_1(X_{k+1},x_{k+1})}
{\left[i_1^{-1}(g)i_2(g)|\ g\in\{{\bf e}_{x_{k+1}}\}\right]}\\
&\cong&\left(\left(\prod\limits_{i=1}^k\pi_1(X_i^*,x_0)\right)*\pi_1(G,x_0)\right)*
\pi_1(X_{k+1},x_{k+1})\\
&\cong&\left(\prod\limits_{i=1}^{k+1}\pi_1(X_i^*,x_0)\right)*\pi_1(G,x_0)\\
&\cong&\left(\prod\limits_{i=1}^m\pi_1(X_i,x_i)\right)*\pi_1(G,x_0),
\end{eqnarray*}
by the induction assumption.\hfill$\Box$

Particularly, for the graph $B_m^T$ or star $S_m^T$ in Fig.$2.2$, we
get the following conclusion.

\vskip 3mm

\no{\bf Corollary $2.5$} \ {\it Let $G$ be the graph $B_m^T$ or star
$S_m^T$. Then}
\begin{eqnarray*}
\pi_1(\mathscr{X}_m\odot
B_m^T,x_0)&\cong&\left(\prod\limits_{i=1}^m\pi_1(X_i^*,x_0)\right)*\pi_1(B_m^T,x_0)\\
&\cong&\left(\prod\limits_{i=1}^m\pi_1(X_i,x_{i-1})\right)*\left<L_i|1\leq
i\leq m\right>,
\end{eqnarray*}
{\it where $L_i$ is the loop of parallel edges $(x_0,x_i)$ in
$B_m^T$ for integers $1\leq i\leq m-1$ and}
$$\pi_1(\mathscr{X}_m\odot
S_m^T,x_0)\cong\prod\limits_{i=1}^m\pi_1(X_i^*,x_0)\cong\prod\limits_{i=1}^m\pi_1(X_i,x_{i-1}).$$

\vskip 3mm

\no{\bf Corollary $2.6$} \ Let $X=\mathscr{X}_m\odot G$ be a
topological space with simply-connected spaces $X_i$ for integers
$1\leq i\leq m$ and $x_0\in X\cap G$. Then we know that
$$\pi_1(X,x_0)\cong\pi_1(G,x_0).$$

\vskip 8mm

\no{\bf \S $3.$ \ A Generalization of Seifert-Van Kampen Theorem}

\vskip 4mm

\no These results and graph $B_m^T$ shown in Section $2$ enables one
to generalize the Seifert-Van Kampen theorem to the case of $U\cap
V$ maybe not arcwise-connected.

\vskip 4mm

\no{\bf Theorem $3.1$} \ {\it Let $X=U\cup V$, $U, V\subset X$ be
open subsets, $X,\ U,\ V$ arcwise-connected and let $C_1,
C_2,\cdots, C_m$ be arcwise-connected components in $U\cap V$ for an
integer $m\geq 1$, $x_{i-1}\in C_i$, $b(x_0,x_{i-1})\subset V$ an
arc $: I\rightarrow X$ with $b(0)=x_0, b(1)=x_{i-1}$ and
$b(x_0,x_{i-1})\cap U=\{x_0,x_{i-1}\}$, $C_i^E=C_i\bigcup
b(x_0,x_{i-1})$ for any integer $i,\ 1\leq i\leq m$, $H$ a group and
there are homomorphisms}
$$
\phi_1^i:\pi_1(U\bigcup b(x_0,x_{i-1}),x_0)\rightarrow H, \ \
\phi_2^i:\pi_1(V,x_0)\rightarrow H$$
{such that}

\c{\includegraphics[bb=10 10 180 120]{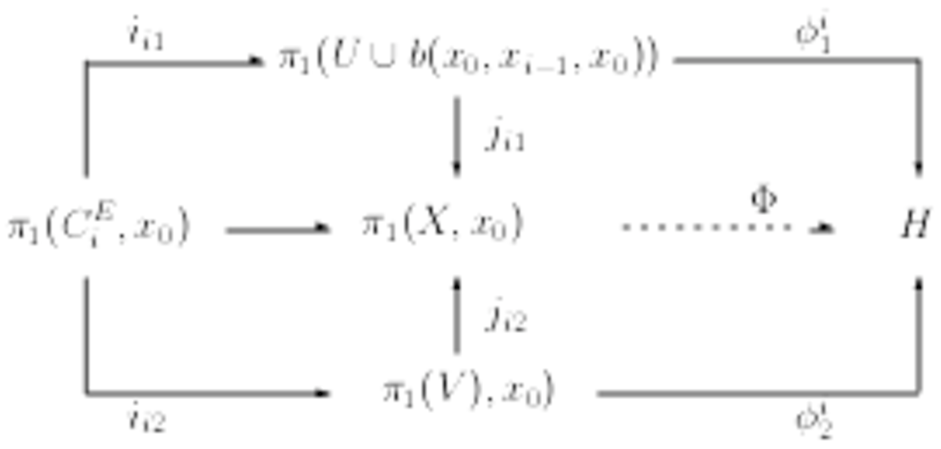}}\vskip 3mm

\no{\it with $\phi_{1}^i\cdot i_{i1}=\phi_{2}^i\cdot i_{i2}$, where
$i_{i1}:\pi_1(C_i^E,x_{0})\rightarrow\pi_1(U\cup
b(x_0,x_{i-1}),x_0)$,
$i_{i2}:\pi_1(C_i^E,x_{0})\rightarrow\pi_1(V,x_0)$ and
$j_{i1}:\pi_1(U\cup b(x_0,x_{i-1},x_0))\rightarrow\pi_1(X,x_0)$,
$j_{i2}:\pi_1(V,x_0))\rightarrow\pi_1(X,x_0)$ are homomorphisms
induced by inclusion mappings, then there exists a unique
homomorphism $\Phi:\ \pi_1(X,x_0)\rightarrow H$ such that $\Phi\cdot
j_{i1}=\phi_{1}^i$ and $\Phi\cdot j_{i2}=\phi_{2}^i$ for integers
$1\leq i\leq m$.} \vskip 3mm

{\it Proof} \ Define $U^{E}=U\bigcup\{\ b(x_0,x_i)\ |\ 1\leq i\leq
m-1\}$. Then we get that $X=U^{E}\cup V$, $U^{E}, V\subset X$ are
still opened with an arcwise-connected intersection $U^{E}\cap
V=\mathscr{X}_m\odot S_m^T$, where $S_m^T$ is a graph formed by arcs
$b(x_0,x_{i-1})$, $1\leq i\leq m$.

Notice that $\mathscr{X}_m\odot Sm^T=\bigcup\limits_{i=1}^mC_i^E$
and $C_i^E\bigcap C_j^E=\{x_0\}$ for $1\leq i,j\leq m,\ i\not=j$.
Therefore, we get that
$$\pi_1(\mathscr{X}_m\odot
S_m^T,x_0)=\bigotimes\limits_{i=1}^m\pi_1(C_i^E,x_0).$$
This fact
enables us knowing that there is a unique $m$-tuple
$(h_1,h_2,\cdots,h_m)$, $h_i\in\pi_1(C_i^E,x_{i-1}),\ 1\leq i\leq m$
such that
$$
\mathscr{I}=\prod\limits_{i=1}^mh_i
$$
for $\forall \mathscr{I}\in\pi_1(\mathscr{X}_m\odot S_m^T,x_0)$.

By definition,
$$i_{i1}:\pi_1(C_i^E,x_{0})\rightarrow\pi_1(U\cap
b(x_0,x_{i-1}),x_0),$$
$$i_{i2}:\pi_1(C_i^E,x_{0})\rightarrow\pi_1(V,x_0)$$
are homomorphisms induced by inclusion mappings. We know that there
are homomorphisms
$$i_1^E:\pi_1(\mathscr{X}_m\odot
S_m^T,x_0)\rightarrow\pi_1(U^E,x_0),$$
$$i_2^E:\pi_1(\mathscr{X}_m\odot
S_m^T,x_0)\rightarrow\pi_1(V,x_0)$$
 with $i_1^E|_{\pi_1(C_i^E,x_0)}=i_{i1}$,
$i_2^E|_{\pi_1(C_i^E,x_0)}=i_{i2}$ for integers $1\leq i\leq m$.

Similarly, because of
$$\pi_1(U^E,x_0)=\bigcup\limits_{i=1}^m\pi_1(U\cup
b(x_0,x_{i-1},x_0))$$
 and
$$j_{i1}:\pi_1(U\cup b(x_0,x_{i-1},x_0))\rightarrow\pi_1(X,x_0),$$
$$j_{i2}:\pi_1(V\rightarrow\pi_1(X,x_0)$$
 being homomorphisms induced by inclusion mappings, there are
homomorphisms
$$j_1^E:\pi_1(U^E,x_0)\rightarrow\pi_1(X,x_0), \ \
j_2^E:\pi_1(V,x_0)\rightarrow\pi_1(X,x_0)$$
 induced by inclusion mappings with $j_1^E|_{\pi_1(U\cup
b(x_0,x_{i-1},x_0))}=j_{i1}$, $j_2^E|_{\pi_1(V,x_0)}=j_{i2}$ for
integers $1\leq i\leq m$ also.

Define $\phi_1^E$ and $\phi_2^E$ by
$$\phi_1^E(\mathscr{I})=\prod\limits_{i=1}^m\phi_1^i(i_{i1}(h_i)), \ \
\phi_2^E(\mathscr{I})=\prod\limits_{i=1}^m\phi_2^i(i_{i2}(h_i)).$$
 Then they are naturally homomorphic extensions of homomorphisms
$\phi_1^i,\ \phi_2^i$ for integers $1\leq i\leq m$. Notice that
$\phi_{1}^i\cdot i_{i1}=\phi_{2}^i\cdot i_{i2}$ for integers $1\leq
i\leq m$, we get that
\begin{eqnarray*}
\phi_1^{E}\cdot i_1^E(\mathscr{I})&=& \phi_1^{E}\cdot i_1^E
\left(\prod\limits_{i=1}^mh_i\right)\\
&=&\prod\limits_{i=1}^m\left(\phi_1^i\cdot i_{i1}(h_i)\right)
=\prod\limits_{i=1}^m\left(\phi_2^i\cdot i_{i2}(h_i)\right)\\
&=&\phi_2^{E}\cdot
i_2^E\left(\prod\limits_{i=1}^mh_i\right)=\phi_2^{E}\cdot
i_2^E(\mathscr{I}),
\end{eqnarray*}
 i.e., the following diagram

\c{\includegraphics[bb=10 10 180 180]{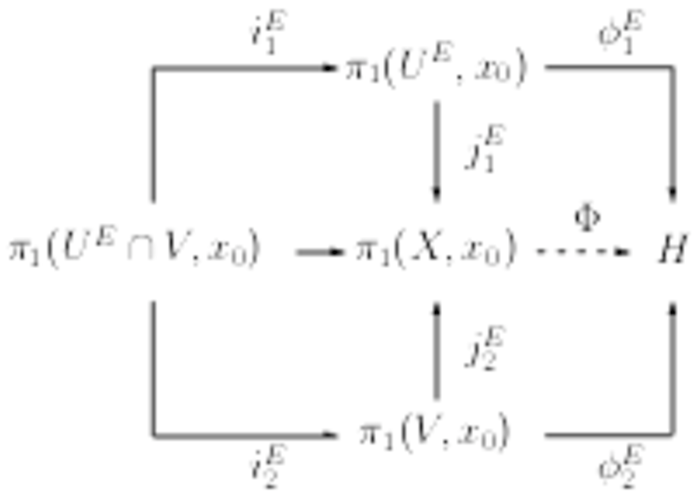}}\vskip 3mm

\no is commutative with $\phi_1^{E}\cdot i_1^E=\phi_2^{E}\cdot
i_2^E$. Applying Theorem $1.1$, we know that there exists a unique
homomorphism $\Phi:\ \pi_1(X,x_0)\rightarrow H$ such that $\Phi\cdot
j_1^E=\phi_1^E$ and $\Phi\cdot j_2^E=\phi_2^E$. Whence,  $\Phi\cdot
j_{i1}=\phi_{1}^i$ and $\Phi\cdot j_{i2}=\phi_{2}^i$ for integers
$1\leq i\leq m$.\hfill$\Box$\vskip 2mm

The following result is a generalization of the classical
Seifert-Van Kampen theorem to the case of maybe non-arcwise
connected.

\vskip 4mm

\no{\bf Theorem $3.2$} \ {\it Let $X$, $U$, $V$, $C_i^E$,
$b(x_0,x_{i-1})$ be arcwise-connected spaces for any integer $i,\
1\leq i\leq m$ as in Theorem $3.1$, $U^{E}=U\bigcup\{\ b(x_0,x_i)\
|\ 1\leq i\leq m-1\}$ and $B_m^T$ a graph formed by arcs
$a(x_0,x_{i-1})$, $b(x_0,x_{i-1})$, $1\leq i\leq m$, where
$a(x_0,x_{i-1})\subset U$ is an arc $: I\rightarrow X$ with
$a(0)=x_0, a(1)=x_{i-1}$ and $a(x_0,x_{i-1})\cap V=\{x_0,x_{i-1}\}$.
Then}

$$\pi_1(X,x_0)\cong\frac{\pi_1(U,x_0)*\pi_1(V,x_0)*\pi_1(B_m^T,x_0)}
{\left[(i_1^{E})^{-1}(g)\cdot i_2(g) |\
g\in\prod\limits_{i=1}^m\pi_1(C_i^E,x_0)\ \right]},$$ {\it where
$i_1^E:\pi_1(U^E\cap V,x_0)\rightarrow\pi_1(U^E,x_0)$ and
$i_2^E:\pi_1(U^E\cap V,x_0)\rightarrow\pi_1(V,x_0)$ are
homomorphisms induced by inclusion mappings.}

\vskip 3mm

{\it Proof} \ Similarly, $X=U^{E}\cup V$, $U^{E}, V\subset X$ are
opened with $U^{E}\cap V=\mathscr{X}_m\odot S_m^T$. By the proof of
Theorem $3.1$ we have known that there are homomorphisms $\phi_1^E$
and $\phi_2^E$ such that $\phi_1^{E}\cdot i_1^E=\phi_2^{E}\cdot
i_2^E$. Applying Theorem $1.2$, we get that
$$\pi_1(X,x_0)\cong\frac{\pi_1(U^E,x_0)*\pi_1(V,x_0)}
{\left[(i_1^{E})^{-1}(\mathscr{I})\cdot i_2^E(\mathscr{I})
|\mathscr{I}\in\pi_1(U^E\cap V,x_0)\right]}.$$

Notice that $U^E\cap V^E=\mathscr{X}_m\odot S_m^T$. We have known
that
$$\pi_1(U^E,x_0)\cong\pi_1(U,x_0)*\pi_1(B_m^T,x_0)$$
 by Corollary $2.3$. As we have shown in the proof of Theorem
$3.1$, an element $\mathscr{I}$ in $\pi_1(\mathscr{X}_m\odot
S_m^T,x_0)$ can be uniquely represented by
$$
\mathscr{I}=\prod\limits_{i=1}^mh_i,
$$
 where $h_i\in\pi_1(C_i^E,x_0),\ 1\leq i\leq m$. We finally get
that
$$\hskip
20mm\pi_1(X,x_0)\cong\frac{\pi_1(U,x_0)*\pi_1(V,x_0)*\pi_1(B_m^T,x_0)}
{\left[(i_1^{E})^{-1}(g)\cdot i_2^E(g) |\
g\in\prod\limits_{i=1}^m\pi_1(C_i^E,x_0)\ \right]}.\hskip 20mm\Box$$

The form of elements in $\pi_1(\mathscr{X}_m\odot S_m^T,x_0)$
appeared in Corollary $2.5$ enables one to obtain another
generalization of classical Seifert-Van Kampen theorem following.

\vskip 4mm

\no{\bf Theorem $3.3$} \ {\it Let $X$, $U$, $V$, $C_1, C_2,\cdots,
C_m$ be arcwise-connected spaces, $b(x_0,x_{i-1})$ arcs for any
integer $i,\ 1\leq i\leq m$ as in Theorem $3.1$, $U^{E}=U\bigcup\{\
b(x_0,x_{i-1})\ |\ 1\leq i\leq m\}$ and $B_m^T$ a graph formed by
arcs $a(x_0,x_{i-1})$, $b(x_0,x_{i-1})$, $1\leq i\leq m$. Then}
$$\pi_1(X,x_0)\cong\frac{\pi_1(U,x_0)*\pi_1(V,x_0)*\pi_1(B_m^T,x_0)}
{\left[(i_1^{E})^{-1}(g)\cdot i_2^E(g) |\
g\in\prod\limits_{i=1}^m\pi_1(C_i,x_{i-1})\right]},$$ {\it where
$i_1^E:\pi_1(U^E\cap V,x_0)\rightarrow\pi_1(U^E,x_0)$ and
$i_2^E:\pi_1(U^E\cap V,x_0)\rightarrow\pi_1(V,x_0)$ are
homomorphisms induced by inclusion mappings.}

\vskip 3mm

{\it Proof} \ Notice that $U^E\cap V=\mathscr{X}_m\odot S^T_m$.
Applying Corollary $2.5$, replacing
$$\pi_1(\mathscr{X}_m\odot S_m^T,x_0)=\left[(i_1^{E})^{-1}(g)\cdot
i_2^E(g) |\ g\in\prod\limits_{i=1}^m\pi_1(C_i^E,x_0)\right]$$
 by
$$\pi_1(\mathscr{X}_m\odot S_m^T,x_0)=\left[(i_1^{E})^{-1}(g)\cdot i_2^E(g) |\
g\in\prod\limits_{i=1}^m\pi_1(C_i,x_{i-1})\right]$$
 in the proof of Theorem $3.2$. We get this
conclusion.\hfill$\Box$\vskip 2mm

Particularly, we get corollaries following by Theorems $3.1$, $3.2$
and $3.3$.

\vskip 4mm

\no{\bf Corollary $3.4$} \ {\it Let $X=U\cup V$, $U, V\subset X$ be
open subsets and $X,\ U,\ V$ and $U\cap V$ arcwise-connected. Then}
$$\pi_1(X,x_0)\cong\frac{\pi_1(U,x_0)*\pi_1(V,x_0)}
{\left[i_1^{-1}(g)\cdot i_2(g)|\ g\in\pi_1(U\cap V,x_0)\right]},$$
{\it where $i_1:\pi_1(U\cap V,x_0)\rightarrow\pi_1(U,x_0)$ and
$i_2:\pi_1(U\cap V,x_0)\rightarrow\pi_1(V,x_0)$ are homomorphisms
induced by inclusion mappings.}

\vskip 3mm

\no{\bf Corollary $3.5$} \ {\it  Let $X$, $U$, $V$, $C_i$,
$a(x_0,x_{i})$, $b(x_0,x_{i})$ for integers $i,\ 1\leq i\leq m$ be
as in Theorem $3.1$. If each $C_i$ is simply-connected, then}
$$\pi_1(X,x_0)\cong\pi_1(U,x_0)*\pi_1(V,x_0)*\pi_1(B_m^T,x_0).$$\vskip
2mm

{\it Proof} \ Notice that $C_1^E, C_2^E,\cdots, C_m^E$ are all
simply-connected by assumption. Applying Theorem $3.3$, we easily
get this conclusion. \hfill$\Box$

\vskip 4mm

\no{\bf Corollary $3.6$} \ {\it  Let $X$, $U$, $V$, $C_i$,
$a(x_0,x_{i})$, $b(x_0,x_{i})$ for integers $i,\ 1\leq i\leq m$ be
as in Theorem $3.1$. If $V$ is simply-connected, then}
$$\pi_1(X,x_0)\cong\frac{\pi_1(U,x_0)*\pi_1(B_m^T,x_0)}
{\left[(i_1^{E})^{-1}(g)\cdot i_2^E(g) |\
g\in\prod\limits_{i=1}^m\pi_1(C_i^E,x_0)\ \right]},$$ {\it where
$i_1^E:\pi_1(U^E\cap V,x_0)\rightarrow\pi_1(U^E,x_0)$ and
$i_2^E:\pi_1(U^E\cap V,x_0)\rightarrow\pi_1(V,x_0)$ are
homomorphisms induced by inclusion mappings.}

\vskip 8mm

\no{\bf \S $4.$ \ Fundamental Groups of Combinatorial Spaces}

\vskip 4mm

\no{\bf 4.1 \ Fundamental groups of combinatorial manifolds}

\vskip 3mm

\no By definition, a combinatorial manifold $\widetilde{M}$ is
arcwise-connected. So we can apply Theorems $3.2$ and $3.3$ to find
its fundamental group $\pi_1(\widetilde{M})$ up to isomorphism in
this section.

\vskip 4mm

\no{\bf Definition $4.1$} \ {\it Let $\widetilde{M}$ be a
combinatorial manifold underlying a graph $G[\widetilde{M}]$. An
edge-induced graph $G^{\theta}[\widetilde{M}]$ is defined by}\vskip
3mm

\hskip 8mm$V(G^{\theta}[\widetilde{M}])=\{x_M,x_{M'},x_1,x_2,\cdots,
x_{\mu(M,M')}|\ for \ \forall (M,M')\in
E(G[\widetilde{M}])\},$\vskip 2mm

\hskip
8mm$E(G^{\theta}[\widetilde{M}])=\{(x_M,x_{M'}),(x_M,x_i),(x_{M'},x_i)|\
1\leq i\leq\mu(M,M')\},$\vskip 2mm

\no{\it where $\mu(M,M')$ is called the edge-index of $(M,M')$ with
$\mu(M,M')+1$ equal to the number of arcwise-connected components in
$M\cap M'$.}\vskip 3mm

By the definition of edge-induced graph, we finally get
$G^{\theta}[\widetilde{M}]$ of a combinatorial manifold
$\widetilde{M}$ if we replace each edge $(M,M')$ in
$G[\widetilde{M}]$ by a subgraph $TB_{\mu(M,M')}^T$ shown in
Fig.$4.1$ with $x_{M}=M$ and $x_{M'}=M'$.

\c{\includegraphics[bb=10 10 180 110]{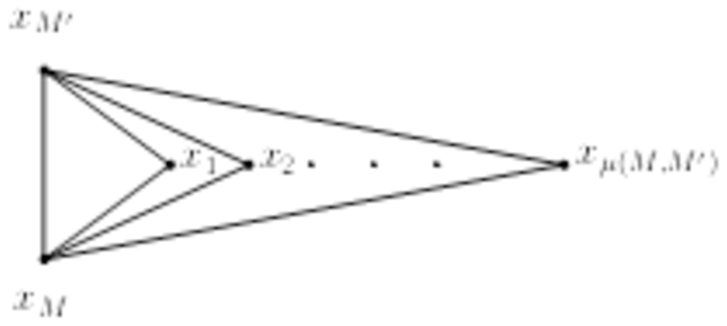}}\vskip 3mm

\c{\bf Fig.$4.1$}\vskip 2mm

Then we have the following result.

\vskip 4mm

\no{\bf Theorem $4.2$} \ {\it Let $\widetilde{M}$ be a finitely
combinatorial manifold. Then}
$$\pi_1(\widetilde{M})\cong\frac{\left(\prod\limits_{M\in
V(G[\widetilde{M}])}\pi_1(M)\right)*\pi_1(G^{\theta}[\widetilde{M}])}
{\left[(i_1^E)^{-1}(g)\cdot i_2^E(g) |\
g\in\prod\limits_{(M_1,M_2)\in E(G[\widetilde{M}])}\pi_1(M_1\bigcap
M_2)\right]},$$ {\it where $i_1^E$ and $i_2^E$ are homomorphisms
induced by inclusion mappings $i_M: \pi_1(M\cap M')\rightarrow
\pi_1(M)$, $i_{M'}:\pi_1(M\cap M')\rightarrow\pi_1(M')$ such as
those shown in the following diagram:

\c{\includegraphics[bb=10 10 180 140]{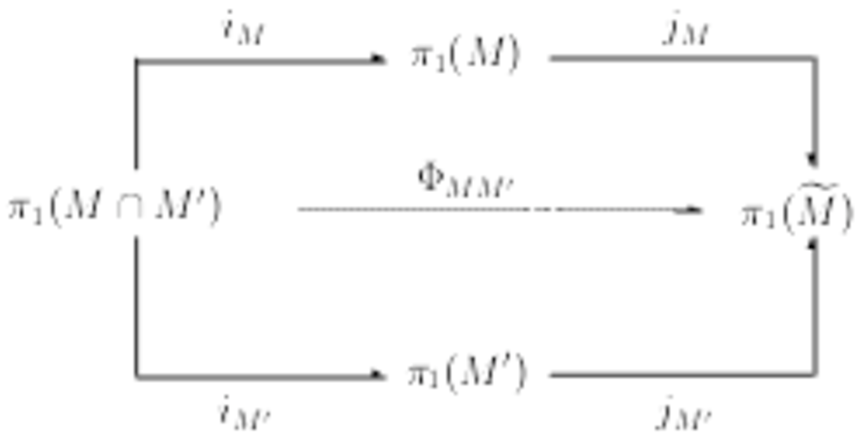}}\vskip 3mm

\no for $\forall(M,M')\in E(G[\widetilde{M}])$.}\vskip 3mm

{\it Proof} \ This result is obvious for $|G[\widetilde{M}]|=1$.
Notice that $G^{\theta}[\widetilde{M}]= B_{\mu(M,M')+1}^T$ if
$V(G[\widetilde{M}])=\{M,\ M'\}$. Whence, it is an immediately
conclusion of Theorem $3.2$ for $|G[\widetilde{M}]|=2$.

Now let $k\geq 3$ be an integer. If this result is true for
$|G[\widetilde{M}]|\leq k$, we prove it hold for
$|G[\widetilde{M}]|=k$. It should be noted that for an
arcwise-connected graph $H$ we can always find a vertex $v\in V(H)$
such that $H-v$ is also arcwise-connected. Otherwise, each vertex
$v$ of $H$ is a cut vertex. There must be $|H|=1$, a contradiction.
Applying this fact to $G[\widetilde{M}]$, we choose a manifold $M\in
V(G[\widetilde{M}])$ such that $\widetilde{M}-M$ is
arcwise-connected, which is also a finitely combinatorial manifold.

Let $U=\widetilde{M}\setminus(M\setminus\widetilde{M})$ and $V=M$.
By definition, they are both opened. Applying Theorem $3.2$, we get
that
$$\pi_1(\widetilde{M})\cong\frac{\pi_1(\widetilde{M}-M)*\pi_1(M)*\pi_1(B_m^T)}
{\left[(i_1^{E})^{-1}(g)\cdot i_2^E(g) |\
g\in\prod\limits_{i=1}^m\pi_1(C_i)\ \right]},$$
 where $C_i$ is an arcwise-connected component in
$M\cap(\widetilde{M}-M)$ and
$$m=\sum\limits_{(M,M')\in E(G[\widetilde{M}])}\mu(M,M').$$
 Notice that
$$\pi_1(B_{m}^T)\cong\prod\limits_{(M,M')\in
E(G[\widetilde{M}]}\pi_1(TB_{\mu(M,M')}).$$
 By the induction assumption, we know that
$$\pi_1(\widetilde{M}-M)\cong\frac{\left(\displaystyle\prod\limits_{M\in
V(G[\widetilde{M}-M])}\pi_1(M)\right)*\pi_1(G^{\theta}[\widetilde{M}-M])}
{\left[(i_1^E)^{-1}(g)\cdot i_2^E(g) |\ \displaystyle
g\in\prod\limits_{(M_1,M_2)\in E(G[\widetilde{M}-M])}\pi_1(M_1\cap
M_2)\right]},$$
 where $i_1^E$ and $i_2^E$ are homomorphisms induced by inclusion
mappings $i_{M_1}: \pi_1(M_1\cap M_2)\rightarrow \pi_1(M_1)$,
$i_{M_2}:\pi_1(M_1\cap M_2)\rightarrow\pi_1(M_2)$ for
$\forall(M_1,M_2)\in E(G[\widetilde{M}-M])$. Therefore, we finally
get that
\begin{eqnarray*}
\pi_1(\widetilde{M})&\cong&\frac{\pi_1(\widetilde{M}-M)*\pi_1(M)*
\pi_1(B_{m}^T)} {\left[(i_1^{E})^{-1}(g)\cdot i_2^E(g) |\
\displaystyle
g\in\prod\limits_{i=1}^m\pi_1(C_i)\ \right]}\\
&\cong& \frac{\frac{\left(\displaystyle\prod\limits_{M\in
V(G[\widetilde{M}-M])}\pi_1(M)\right)\displaystyle*\pi_1(G^{\theta}[\widetilde{M}-M])}
{\left[\displaystyle(i_1^E)^{-1}(g)\cdot i_2^E(g) |\
g\in\prod\limits_{(M_1,M_2)\in E(G[\widetilde{M}-M])}\pi_1(M_1\cap
M_2)\right]}} {\left[\displaystyle(i_1^{E})^{-1}(g)\cdot i_2(g) |\
g\in\prod\limits_{i=1}^m\pi_1(C_i)\
\right]}\\
&*&\frac{\pi_1(M)*\displaystyle\prod\limits_{(M,M')\in
E(G[\widetilde{M}]}\pi_1(TB_{\mu(M,M')})}{\left[\displaystyle(i_1^{E})^{-1}(g)\cdot
i_2(g) |\ g\in\prod\limits_{i=1}^m\pi_1(C_i)\ \right]}\\
&\cong&\frac{\left(\displaystyle\prod\limits_{M\in
V(G[\widetilde{M}])}\pi_1(M)\right)*\pi_1(G^{\theta}[\widetilde{M}])}
{\left[(i_1^E)^{-1}(g)\cdot i_2^E(g) |\ \displaystyle
g\in\prod\limits_{(M_1,M_2)\in E(G[\widetilde{M}])}\pi_1(M_1\bigcap
M_2)\right]}
\end{eqnarray*}
 by facts
$$\left(\mathscr{G}/\mathscr{H}\right)*H\cong\mathscr{G}*H/\mathscr{H}$$
 for groups $\mathscr{G,\ H}$, $G$ and
$$G^{\theta}[\widetilde{M}]=G^{\theta}[\widetilde{M}-M]\bigcup\limits_{(M,M')\in
E(G[\widetilde{M}]}TB_{\mu(M,M')},$$
$$\pi_1(G^{\theta}[\widetilde{M}])=\pi_1(G^{\theta}[\widetilde{M}-M])*\prod\limits_{(M,M')\in
E(G[\widetilde{M}]}\pi_1(TB_{\mu(M,M')}),
$$
$$\prod\limits_{M\in
V(G[\widetilde{M}])}\pi_1(M)=\left(\prod\limits_{M\in
V(G[\widetilde{M}-M])}\pi_1(M)\right)*\pi_1(M),$$
 where $i_1^E$ and $i_2^E$ are homomorphisms induced by inclusion
mappings $i_M: \pi_1(M\cap M')\rightarrow \pi_1(M)$,
$i_{M'}:\pi_1(M\cap M')\rightarrow\pi_1(M')$ for $\forall(M,M')\in
E(G[\widetilde{M}])$. This completes the proof.\hfill$\Box$

Applying Corollary $3.5$, we get a result known in [8] by noted that
$G^{\theta}[\widetilde{M}]=G[\widetilde{M}]$ if
$\forall(M_1,M_2)\in E(G^L[\widetilde{M}])$, $M_1\cap M_2$ is simply
connected.

\vskip 4mm

\no{\bf Corollary $4.3$}([8]) \ {\it Let $\widetilde{M}$ be a
finitely combinatorial manifold. If for $\forall(M_1,M_2)\in
E(G^L[\widetilde{M}])$, $M_1\cap M_2$ is simply connected, then }
$$
\pi_1(\widetilde{M})\cong \left(\bigoplus\limits_{M\in
V(G[\widetilde{M}])}\pi_1(M)\right)\bigoplus
\pi_1(G[\widetilde{M}]).
$$

\vskip 4mm

\no{\bf 4.2 \ Fundamental groups of manifolds}

\vskip 3mm

\no Notice that $\pi_1({\bf R}^n)=identity$ for any integer $n\geq
1$. If we choose $M\in V(G[\widetilde{M}])$ to be a chart
$(U_{\lambda},\varphi_{\lambda})$ with $\varphi_{\lambda}:
U_{\lambda}\rightarrow{\bf R}^n$ for $\lambda\in\Lambda$ in Theorem
$4.2$, i.e., an $n$-manifold, we get the fundamental group of
$n$-manifold following.

\vskip 4mm

\no{\bf Theorem $4.4$} \ {\it Let $M$ be a compact $n$-manifold with
charts $\{(U_{\lambda},\varphi_{\lambda})|\
\varphi_{\lambda}:U_{\lambda}\rightarrow{\bf R}^n,
\lambda\in\Lambda)\}$. Then}
$$\pi_1(M)\cong\frac{\pi_1(G^{\theta}[M])}
{\left[(i_1^E)^{-1}(g)\cdot i_2^E(g) |\ g\in\prod\limits_{(U_{\mu},
U_{\nu})\in E(G[M])}\pi_1(U_{\mu}\cap U_{\nu})\right]},$$ {\it where
$i_1^E$ and $i_2^E$ are homomorphisms induced by inclusion mappings
$i_{U_{\mu}}: \pi_1(U_{\mu}\cap U_{\nu})\rightarrow \pi_1(U_{\mu})$,
$i_{U_{\nu}}:\pi_1(U_{\mu}\cap U_{\nu})\rightarrow\pi_1(U_{\nu})$,
$\mu,\nu\in\Lambda$.}

\vskip 4mm

\no{\bf Corollary $4.5$} \ {\it Let $M$ be a simply connected
manifold with charts $\{(U_{\lambda},\varphi_{\lambda})|\
\varphi_{\lambda}:U_{\lambda}\rightarrow{\bf R}^n,
\lambda\in\Lambda)\}$, where $|\Lambda|<+\infty$. Then
$G^{\theta}[M]=G[M]$ is a tree.}

\vskip 3mm

Particularly, if $U_{\mu}\cap U_{\nu}$ is simply connected for
$\forall\mu,\nu\in\Lambda$, then we obtain an interesting result
following.

\vskip 4mm

\no{\bf Corollary $4.6$} \ {\it Let $M$ be a compact $n$-manifold
with charts $\{(U_{\lambda},\varphi_{\lambda})|\
\varphi_{\lambda}:U_{\lambda}\rightarrow{\bf R}^n,
\lambda\in\Lambda)\}$. If $U_{\mu}\cap U_{\nu}$ is simply connected
for $\forall\mu,\nu\in\Lambda$, then}
$$\pi_1(M)\cong\pi_1(G[M]).$$

\vskip 3mm

Therefore, by Theorem $4.4$ we know that the fundamental group of a
manifold $M$ is a subgroup of that of its edge-induced graph
$G^{\theta}[M]$. Particularly, if $G^{\theta}[M]=G[\widetilde{M}]$,
i.e., $U_{\mu}\cap U_{\nu}$ is simply connected for
$\forall\mu,\nu\in\Lambda$, then it is nothing but the fundamental
group of $G[\widetilde{M}]$.

\vskip 8mm

\no{\bf References}\vskip 3mm

\re{[1]}Munkres J.R., {\it Topology} (2nd edition), Prentice Hall,
Inc, 2000.

\re{[2]}W.S.Massey, {\it Algebraic Topology: An Introduction},
Springer-Verlag, New York, etc.(1977).

\re{[3]}John M.Lee, {\it Introduction to Topological Manifolds},
Springer-Verlag New York, Inc., 2000.

\re{[4]}J.L.Gross and T.W.Tucker, {\it Topological Graph Theory},
John Wiley \& Sons, 1987.

\re{[5]}L.F.Mao, {\it Automorphism Groups of Maps, Surfaces and
Smarandache Geometries}, American Research Press, 2005.

\re{[6]}L.F.Mao, Geometrical theory on combinatorial manifolds, {\it
JP J.Geometry and Topology}, Vol.7, No.1(2007),65-114.

\re{[7]}L.F.Mao, Combinatorial fields - an introduction, {\it
International J.Math. Combin.} Vol.3 (2009), 01-22.

\re{[8]}L.F.Mao, {\it Combinatorial Geometry with Applications to
Field Theory}, InfoQuest, USA, 2009.

\re{[9]}L.F.Mao, {\it Smarandache Multi-Space Theory}, Hexis,
Phoenix, USA 2006.

\re{[10]}Smarandache F. Mixed noneuclidean geometries. arXiv:
math/0010119, 10/2000.
\end{document}